\documentclass[reqno]{amsart}
\usepackage{amsmath,amssymb,enumerate}

\theoremstyle{plain}
\newtheorem{theorem}{Theorem}
\newtheorem*{theorem*}{Theorem}
\newtheorem{proposition}[theorem]{Proposition}
\newtheorem*{proposition*}{Proposition}

\newtheorem{lemmma*}{Lemma}

\theoremstyle{definition}

\newtheorem*{definition*}{Definition}

\newtheorem*{remark*}{Remark}
\newcommand{\Irr}{\mathrm{Irr}}
\newcommand{\Cl}{\mathrm{Cl}}
\begin{document}
\title{Short character values on large conjugacy~classes}
\author[A.\ R.\ Miller]{Alexander Rossi Miller}
\address{Center for Communications Research, Princeton}
\email{a.miller@ccr-princeton.org}
\begin{abstract}
  We provide an example of a finite group with a conjugacy class of
  average size on which fewer than half of the irreducible
  characters are either zero or a root of unity.
\end{abstract}
\maketitle
\thispagestyle{empty}
For any irreducible character $\chi$ of a finite group $G$, let
\[
  \theta(\chi)=\frac{|\{g\in G : \chi(g)\ \text{is either zero or a root of unity}\}|}{|G|}.
\]
It is a classic result of J.~G.~Thompson that, for any finite group $G$, 
\[
  \theta(\chi)>\tfrac{1}{3},\quad \chi\in \Irr(G).
\]
In a recent paper~\cite{MillerA1} the author conjectured that, in fact, for any finite group $G$,
\begin{equation}\label{Main Conjecture}
  \theta(\chi)\geq \tfrac{1}{2},\quad \chi\in\Irr(G).
\end{equation}
There are no open interval gaps beyond $\frac{1}{2}$ in the sense \cite{MillerA3} that for any
$l\in[\frac{1}{2},1]$ and any $\epsilon>0$, there exists a finite group $G$ and an irreducible character
$\chi$ of $G$ such that $|\theta(\chi)-l|<\epsilon$. 
The conjectured bound \eqref{Main Conjecture} was established in \cite{MillerA1} for many groups,
including all finite nilpotent groups, all sporadic simple groups, all simple groups of order $<10^9$,
and various infinite families of simple groups.
A.~Moret\'o and G.~Navarro \cite{MoretoNavarro} recently established \eqref{Main Conjecture} for all finite groups that admit a~Sylow~series. 

The bound in \eqref{Main Conjecture} has now also been verified for a large number of small groups.

\begin{proposition}
  If $G$ is a group of order $\leq 2000$, then
  $\theta(\chi)\geq \frac{1}{2}$ for all $\chi\in\Irr(G)$.
\end{proposition}

We remark that computing the character tables to verify
\eqref{Main Conjecture} for the $408641062$ groups of order $1536$
with a standard computer algebra package would take a $3$GHz~CPU several years,
but by inspection we find that $10494213$ of these groups of order $1536$ are nilpotent and
$398050412$ of the non-nilpotent groups of order $1536$ admit a Sylow series,
which leaves only $96437$ groups of order $1536$ to check.
These are the groups $\texttt{SmallGroup(1536,i)}$ with ${408 544 626\leq i\leq 408 641 062}$,
and for each of these groups the character table was computed in order to verify~\eqref{Main Conjecture}.

There is also an analogue of Thompson's result for columns of character tables. 
For any finite group $G$, and any element $g\in G$, let
\[
  \theta'(g)=\frac{|\{\chi\in \Irr(G) : \chi(g)\ \text{is either zero or a root of unity}\}|}{|\Irr(G)|}
\]
and let
$L'(G)=\{x\in G : |\Cl(G)|\geq|C_G(x)|\}$.
Then, following Thompson, P.~X.~Gallagher \cite{Gallagher} proved that
\[
  \theta'(g)>\tfrac{1}{3},\quad g\in  L'(G).
\]
It was conjectured in \cite{MillerA1} that
$\theta'(g)\geq \tfrac{1}{2}$ for $g\in L'(G)$.
This turns out to be~too~strong.

\begin{proposition}
  There is a group $G$ of order $1960$ with $\theta'(g)=\frac{16}{35}$ for some~${g\in L'(G)}$.
\end{proposition}

To find this counterexample, we carried out an exhaustive search of the
groups of order $\leq 2000$, excluding those of order $1536$, and found exactly one counterexample. 
This counterexample $X$ is the group \texttt{SmallGroup(1960,122)}. 
It is not nilpotent, but it is monomial and it does admit a Sylow series.\footnote{The group $X$ also has the property that, for each $\chi\in\Irr(X)$, $|\chi(x)|=1$ if and only if~${\chi(1)=1}$. Call any group with this property a \emph{circle-avoiding group}. It was shown in \cite{MillerA1} that nilpotent groups are circle-avoiding. Many circle-avoiding groups are monomial, but this is not always the case. For example, \texttt{SmallGroup(1500,36)} is a circle-avoiding group that is not monomial, but it is solvable.} 
For this group $X$, we find that $\theta'(x)<\frac{1}{2}$ for exactly twelve classes $x^X\subset L'(X)$,
each of which satisfies
$\theta'(x)=\frac{16}{35}$ and ${|C_X(x)|=|\Cl(X)|}$. 
This suggests that perhaps the original conjecture requires only minor modification. 
Is it true that, for any finite group~$G$, 
\begin{equation}\label{Modified conjecture}
  \theta'(g)\geq \tfrac{1}{2},\quad g\in L(G),
\end{equation}
where
$L(G)=\{g\in G : |\Cl(G)|>|C_G(g)|\}$
is the set of all elements belonging to strictly larger-than-average classes?

\end{document}